\documentclass[11pt]{article}
\usepackage{graphicx}
\usepackage{amsfonts}
  \def\fH{{\cal H}}

\usepackage{theorem}
\newtheorem{Lem}{Lemma}[section]

\newtheorem{The}[Lem]{Theorem}
\newtheorem{Prop}[Lem]{Proposition}
\newtheorem{Cor}[Lem]{Corollary}

\newtheorem{Rem}[Lem]{Remark}

\newcommand{\qed}{\hbox{\rule{6pt}{6pt}}}

\begin{document}
\title{Operator inequalities among arithmetic mean, geometric mean and harmonic mean}
\author{Shigeru Furuichi\footnote{E-mail:furuichi@chs.nihon-u.ac.jp}\\
{\small Department of Information Science,}\\
{\small College of Humanities and Sciences, Nihon University,}\\
{\small 3-25-40, Sakurajyousui, Setagaya-ku, Tokyo, 156-8550, Japan}}
\date{}
\maketitle
{\bf Abstract.} 
We give an upper bound for the weighted geometric mean using the weighted arithmetic mean and the weighted harmonic mean. We also give a lower bound for 
the weighted geometric mean. These inequalities are proven for two invertible positive operators.
\vspace{3mm}


\vspace{3mm}
\vspace{3mm}

\section{Introduction}
Let $\fH$ be a complex Hilbert space. We represent the set of all bounded operators on $\fH$ by $B(\fH)$.
If $A\in B(\fH)$ satisfies $A^*=A$, then $A$ is called a self-adjoint operator.
If a self-adjoint operator $A$ satisfies $\langle x \vert A \vert x\rangle \geq 0$ for any $\vert x \rangle \in \fH$, then $A$ is called a positive operator.
For two self-adjoint operators $A$ and $B$, $A\geq B$ means $A-B\geq 0$.
The notation $A >0$ means $A$ is an invertible positive operator.

It is well-known that we have the following Young inequalities for invertible positive operators $A$ and $B$:
\begin{equation} \label{orig_Young_ineq}
(1-\nu) A+ \nu B \geq A\#_{\nu}B \geq \left\{(1-\nu) A^{-1}+ \nu B^{-1} \right\}^{-1},
\end{equation}
where $A\#_{\nu}B \equiv A^{1/2}(A^{-1/2}BA^{-1/2})^{\nu}A^{1/2}$  
represents the geometric mean for two positive operators $A$ and $B$ and a  weighted parameter $\nu\in[0,1]$ \cite{KA}. (In this paper, we  use the notation $A \# B$ instead of $A \#_{1/2} B$ for the simplicity.) 
$(1-\nu) A+ \nu B$ and $\left\{(1-\nu) A^{-1}+ \nu B^{-1} \right\}^{-1}$ are called weighted arithmetic mean and harmonic mean for two positive operators, respectively.
The simplified and elegant proof for the inequalities (\ref{orig_Young_ineq}) 
was given in \cite{FY}. 
Recently, refinements of the inequalities (\ref{orig_Young_ineq}) were given in our papers \cite{Furu1,Furu2}. 
It is also notable that improvements of \cite{Furu2} have been given in the
paper \cite{ZSF}. And further improvements have been given in quite recent papers \cite{KLP} and \cite{KKLP}. 
In this short note, we consider the relations among operator means for two positive operators. 

We start from the following proposition.
\begin{Prop} \label{start_proposition}
Let $A, B$ be invertible positive operators and $r$ be a real number. Then we have the following inequalities.
\begin{itemize}
\item[(i)] If $r \geq 2$, then $rA\# B + \left( {1 - r} \right)\frac{{A + B}}{2} \le {\left( {\frac{{{A^{ - 1}} + {B^{ - 1}}}}{2}} \right)^{ - 1}}.$
\item[(ii)] If $r \leq 1$, then $rA\# B + \left( {1 - r} \right)\frac{{A + B}}{2} \ge {\left( {\frac{{{A^{ - 1}} + {B^{ - 1}}}}{2}} \right)^{ - 1}}.$
\end{itemize}
\end{Prop}

{\it Proof :}
In general, by using the notion of the representing function $f_m(x)=1 m x$ for operator mean $m$, 
it is well-known \cite{KA} that $f_m(x) \leq f_n(x)$ holds for $x>0$ if and only if $A m B \leq A n B$ holds for all positive operators $A$ and $B$. 
Thus we can prove this proposition from the following scalar inequalities for $t>0$.
\begin{itemize}
\item[(i)] $r\sqrt t  + \left( {1 - r} \right)\frac{{t + 1}}{2} \le \frac{{2t}}{{t + 1}},
\quad \left( {r \ge 2} \right).$
\item[(ii)] $r\sqrt t  + \left( {1 - r} \right)\frac{{t + 1}}{2} \ge \frac{{2t}}{{t + 1}},
\quad \left( {r \le 1} \right).$
 \end{itemize}
Actually (i) above can be proven in the following way.
We set $f_r(t) \equiv \frac{2t}{t+1}-r \sqrt{t} -(1-r) \frac{t+1}{2}$, then
$\frac{df_r(t)}{dr} = -\sqrt{t} + \frac{t+1}{2} \geq 0$ implies $f_r(t) \geq f_2(t)$ for $r \geq 2$. From the relation $\frac{2t}{t+1} + \frac{t+1}{2} \geq  2\sqrt{t}$, we have $f_2(t) \geq 0$.
We also give the proof for (ii) above.
We set $g_r(t) \equiv r\sqrt{t} +(1-r)\frac{t+1}{2} - \frac{2t}{t+1}$,
then $\frac{dg_r(t)}{dr} = \sqrt{t}-\frac{t+1}{2} \leq 0$ implies $g_r(t) \geq g_1(t)$
for $r \leq 1$. From the relation $\frac{2t}{t+1}\leq \sqrt{t}$, we have $g_1(t) \geq 0$.

\hfill \qed

\begin{Rem}
We have  counter-examples of both inequalities (i) and (ii) in Proposition \ref{start_proposition} for $1 < r < 2$.
For example, we take $r=1.5$. Then we have the following computations.
$\frac{2t}{t+1} -r \sqrt{t} -(1-r)\frac{t+1}{2} \simeq 0.122302$ when $t=0.01$
and $\frac{2t}{t+1} -r \sqrt{t} -(1-r)\frac{t+1}{2} \simeq -0.037987$ when $t=2$.
\end{Rem}


\section{Main results}
Proposition \ref{start_proposition} can be generalized by means of weighted parameter $\nu \in [0,1]$, as the second inequality in  (\ref{main_theorem_ineq00}) below.
\begin{The} \label{main_theorem00}
If (i) $0 \leq \nu \leq 1/2$ and $0 < A \leq B$ or (ii)
$1/2 \leq \nu \leq1$ and $0 < B \leq A$, then the following inequalities hold
\begin{equation} \label{main_theorem_ineq00}
A\# B + \left( {\nu  - \frac{1}{2}} \right)\left( {B - A} \right) \le A{\# _\nu }B \le \frac{1}{2}\left\{ {\left( {1 - \nu } \right)A + \nu B} \right\} + \frac{1}{2}{\left\{ {\left( {1 - \nu } \right){A^{ - 1}} + \nu {B^{ - 1}}} \right\}^{ - 1}}.
\end{equation}
\end{The}

\begin{Rem}
Under the same conditions as in Theorem \ref{main_theorem00}, we have $A\# B \ge A{\# _\nu }B.$
\end{Rem}
In order to prove Theorem \ref{main_theorem00}, we firstly prove the corresponding scalar inequalities, as it was similarly done in  Proposition \ref{start_proposition}.
 
\begin{Lem} \label{sec2_lemma01}
If (i) $0 \leq  \nu \leq 1/2$ and $t \geq 1$ or (ii)
$1/2 \leq  \nu \leq 1$ and $0 < t \leq  1$, then the following inequalities hold
\begin{equation} \label{key_lemma_ineq01}
2\sqrt t  + \left( {2\nu  - 1} \right)\left( {t - 1} \right) \le 2{t^\nu } \le \left( {1 - \nu } \right) + \nu t + {\left\{ {\left( {1 - \nu } \right) + \frac{\nu }{t}} \right\}^{ - 1}}.
\end{equation}
\end{Lem}

{\it Proof:}
It is trivial for the case $t=1$. For the cases $\nu =0$, $1/2$ or $1$, the inequalities (\ref{key_lemma_ineq01}) hold. So we assume $t \ne 1$ and $\nu \ne 0,1/2,1$.
We firstly prove the first inequality of the inequalities (\ref{key_lemma_ineq01}),
under the condition (i) $0 < \nu < 1/2$ and $t > 1$ or (ii)
$1/2 < \nu < 1$ and $0 <  t < 1$.
Here we put ${f_\nu }\left( t \right) \equiv {t^\nu } - \sqrt t  - \left( {\nu  - \frac{1}{2}} \right)\left( {t - 1} \right).$ Then we have 
$f{'_\nu }\left( t \right) = \nu {t^{\nu  - 1}} - \frac{1}{2}\frac{1}{{\sqrt t }} - \left( {\nu  - \frac{1}{2}} \right)$ and $f{'_\nu }\left( 1 \right) = 0$. We also have
$f'{'_\nu }\left( t \right) =  - \nu \left( {1 - \nu } \right){t^{\nu  - 2}} + \frac{1}{4}{t^{ - 3/2}}$. Thus we have 
$f'{'_\nu }\left( t \right) = 0 \Leftrightarrow t = {t_\nu } \equiv {\left\{ {4\nu \left( {1 - \nu } \right)} \right\}^{\frac{2}{{1 - 2\nu }}}}$.
We find $t_{\nu} < 1 $ in the case $0 < \nu <1/2$ and $t>1$.
Then we find $f'{'_\nu }\left( t \right) \ge 0$ for $t>1 (> t_{\nu})$.
So $f{'_\nu }\left( t \right)$ is monotone increasing for $t > 1$ and we have $f{'_\nu }\left( 1 \right) = 0.$ Thus we find $f{'_\nu }\left( t \right) \ge 0$ for $t > 1$. So $f{_\nu }\left( t \right)$ is monotone increasing for $t > 1$.  Therefore
we have ${f_\nu }\left( t \right) \ge {f_\nu }\left( 1 \right) = 0.$
We also find $t_{\nu} >1 $ in the case $1/2 < \nu <1$ and $0<t<1$.
Then we find $f'{'_\nu }\left( t \right) \ge 0$ for $0 < t < 1 (< t_{\nu})$.
So $f{'_\nu }\left( t \right)$ is monotone increasing for $0 < t < 1$ and we have $f{'_\nu }\left( 1 \right) = 0.$ Thus we find $f{'_\nu }\left( t \right) \le 0$ for $0<t<1$. So  $f{_\nu }\left( t \right)$ is monotone decreasing for $0 < t < 1$.  Therefore we have ${f_\nu }\left( t \right) \ge {f_\nu }\left( 1 \right) = 0.$
Thus the proof for the first inequality of the inequalities (\ref{key_lemma_ineq01}) is done.

We prove the second inequality of the inequalities (\ref{key_lemma_ineq01}). 
We put ${g_\nu }\left( t \right) \equiv \left( {1 - \nu } \right) + \nu t + \frac{1}{{1 - \nu  + \frac{\nu }{t}}} - 2{t^\nu }$. Then we have
${g_\nu }\left( t \right) = \left( {1 - \nu } \right) + \nu t + \frac{t}{{\left( {1 - \nu } \right)t + \nu }} - 2{t^\nu } \ge 2\sqrt {\frac{{\left\{ {\left( {1 - \nu } \right) + \nu t} \right\}t}}{{\left( {1 - \nu } \right)t + \nu }}}  - 2{t^\nu }$.
Since ${g_\nu }\left( t \right) \ge 0$ is equivalent to $\frac{{\left( {1 - \nu } \right) + \nu t}}{{\left( {1 - \nu } \right)t + \nu }} \ge {t^{2\nu  - 1}}$, we put again
${h_\nu }\left( t \right) \equiv \left( {1 - \nu } \right) + \nu t - \left\{ {\left( {1 - \nu } \right)t + \nu } \right\}{t^{2\nu  - 1}}$.
Then we prove $h_{\nu}(t) >0$ under the condition (i) $0 < \nu <1/2$ and $t>1$ or (ii) $1/2 < \nu <1$ and $0<t<1$.
By the elementary calculations, we have
${h_\nu }'\left( t \right) = \nu  - 2\nu \left( {1 - \nu } \right){t^{2\nu  - 1}} - \nu \left( {2\nu  - 1} \right){t^{2\nu  - 2}}$, ${h_\nu }'\left( 1 \right) = 0$ and ${h_\nu }''\left( t \right) =  - 2\nu \left( {1 - \nu } \right)\left( {2\nu  - 1} \right){t^{2\nu  - 3}}\left( {t - 1} \right).$ Then we find ${h_\nu }''\left( t \right) = 0 \Leftrightarrow t = 1$.
In the case $t > 1$, we have ${h_\nu }''\left( t \right) \ge 0.$
So ${h_\nu }'\left( t \right)$ is monotone increasing for $ t > 1$  and we have ${h_\nu }'\left( 1 \right) = 0$. Thus we have ${h_\nu }'\left( t \right) \ge 0$ for $t > 1$. So ${h_\nu }\left( t \right)$ is monotone increasing for $ t > 1$. Thus we have ${h_\nu }\left( t \right) \ge {h_\nu }\left( 1 \right) = 0$.
In the case $0 < t <1$, we also have ${h_\nu }''\left( t \right) \ge 0.$
So ${h_\nu }'\left( t \right)$ is monotone increasing for $ 0 < t < 1$  and we have ${h_\nu }'\left( 1 \right) = 0$. Thus we have ${h_\nu }'\left( t \right) \leq 0$ for 
$0 < t < 1$. So ${h_\nu }\left( t \right)$ is monotone decreasing for $ 0 < t < 1$. Thus we have ${h_\nu }\left( t \right) \ge {h_\nu }\left( 1 \right) = 0$.
Thus the proof for the second inequality of the inequalities (\ref{key_lemma_ineq01}) is done.

\hfill \qed

\begin{Lem} \label{sec2_lemma02}
Let $r \in \mathbb{R}$.
Then the function ${k_{r,\nu }}\left( t \right) \equiv r{t^\nu } + \left( {1 - r} \right)\left\{ {\left( {1 - \nu } \right) + \nu t} \right\},\left( {0 \leq \nu   \leq 1,   t > 0} \right)$ is monotone decreasing with respect to $r$. Therefore,
${k_{r,\nu }}\left( t \right) \le {k_{2,\nu }}\left( t \right)$ for $r \ge 2$
and ${k_{r,\nu }}\left( t \right) \ge {k_{1,\nu }}\left( t \right)$ for $r \le 1$.
\end{Lem}
{\it Proof :}
The proof is done by $\frac{{\partial {k_{r,\nu }}\left( t \right)}}{{\partial r}} = {t^\nu } - \left\{ {\left( {1 - \nu } \right) + \nu t} \right\} \le 0$, for  $\nu \in [0,1]$ and  $t>0$.

\hfill \qed

Lemma \ref{sec2_lemma02} provides the following results.
\begin{Lem} \label{sec2_lemma03}
Let $r \geq 2$. If (i) $0  \leq \nu  \leq1/2$ and $t \geq 1$ or (ii)
$1/2  \leq \nu  \leq1$ and $0<t \leq 1$, then 
$$r{t^\nu } + \left( {1 - r} \right)\left\{ {\left( {1 - \nu } \right) + \nu t} \right\} \le {\left\{ {\left( {1 - \nu } \right) + \frac{\nu }{t}} \right\}^{ - 1}}.$$
\end{Lem}
{\it Proof :}
The proof follows directly from Lemma \ref{sec2_lemma01} and  Lemma \ref{sec2_lemma02}.

\hfill \qed

\begin{Lem}\label{sec2_lemma04}
Let $r \leq  1$.  For $0 < \nu   \leq 1$ and $t > 0$, we have
\[r{t^\nu } + \left( {1 - r} \right)\left\{ {\left( {1 - \nu } \right) + \nu t} \right\} \ge {\left\{ {\left( {1 - \nu } \right) + \frac{\nu }{t}} \right\}^{ - 1}}.\]
\end{Lem}
{\it Proof:}
For $r \leq 1$, it follows from Lemma \ref{sec2_lemma02} that $r{t^\nu } + \left( {1 - r} \right)\left\{ {\left( {1 - \nu } \right) + \nu t} \right\} \ge {t^\nu }$.
Since we have ${t^\nu } \ge {\left\{ {\left( {1 - \nu } \right) + \frac{\nu }{t}} \right\}^{ - 1}}$, the proof is done.

\hfill \qed

Finally we have the following corollary. 
\begin{Cor}
Let $r \geq 2$. If  (i) $0 < \nu  \leq 1/2$ and $0 < A \leq B$ or (ii) $1/2  \leq \nu  \leq 1$ and $0 < B \leq A$, then 
\[rA{\# _\nu }B + \left( {1 - r} \right)\left\{ {\left( {1 - \nu } \right)A + \nu B} \right\} \le {\left\{ {\left( {1 - \nu } \right){A^{ - 1}} + \nu {B^{ - 1}}} \right\}^{ - 1}}.\]

Let $r \leq 1$. For $0 < \nu <1$ and $t >0$, we have
\[rA{\# _\nu }B + \left( {1 - r} \right)\left\{ {\left( {1 - \nu } \right)A + \nu B} \right\} \ge {\left\{ {\left( {1 - \nu } \right){A^{ - 1}} + \nu {B^{ - 1}}} \right\}^{ - 1}}.\]
\end{Cor}
{\it Proof:}
The proof can be done applying Lemma \ref{sec2_lemma03}, Lemma \ref{sec2_lemma04} and Theorem \ref{main_theorem00}.

\hfill \qed


\section*{Acknowledgement}
The author was partially supported by JSPS KAKENHI Grant Number 24540146.


\begin{thebibliography}{99}
\bibitem{KA} F.Kubo and T.Ando, Means of positive operators, Math. Ann.,Vol.264(1980),pp.205-224.
\bibitem{FY} T.Furuta and M.Yanagida, Generalized means and convexity of inversion for positive operators, Amer.Math.Monthly, Vol.105 (1998),pp.258-259.
\bibitem{Furu1} S.Furuichi, On refined Young inequalities and reverse inequalities, J. Math. Ineq., Vol.5(2011), pp.21-31.
\bibitem{Furu2} S.Furuichi, Refined Young inequalities with Specht's ratio, J. Egypt. Math. Soc., Vol.20(2012), pp.46-49.
\bibitem{ZSF} H.Zuo, G.Shi and M.Fujii, Refined Young inequality with Kantrovich constant, J. Math. Ineq., Vol.5(2011), pp.551-556.
\bibitem{KLP} M.Krni\'c, N.Lovci\v{c}evi\'c and J.Pe\v{c}ari\'c, Jensen's operator and applications to mean inequalities for operators in Hilbert space, Bull. Malays. Math. Sci. Soc., Vol.35 (2012), pp.1-14. 
\bibitem{KKLP} F.Kittanech, M.Krni\'c, N.Lovci\v{c}evi\'c and J.Pe\v{c}ari\'c,
Improved arithmetic-geometric and Heinz means inequalities for Hilbert space operators, Publ. Math. Debrecen, Vol.80(2012), pp.465-478.
\end{thebibliography}
\end{document}